\documentclass[12pt]{article}
\usepackage{amssymb}
\usepackage{amsfonts}
\usepackage{amsmath}
\usepackage{graphicx}

\providecommand{\U}[1]{\protect\rule{.1in}{.1in}}
\newtheorem{theorem}{Theorem}[section]

\newtheorem{lemma}[theorem]{Lemma}

\newtheorem{conclusion}[theorem]{Conclusion}
\input{tcilatex}
\begin{document}

\title{A fast method of reconstruction for X-ray phase contrast imaging with
arbitrary Fresnel number.}
\author{Victor Palamodov \and Tel Aiv University,}
\maketitle

\textbf{Abstract. }New lensless diffractive X-ray technic for micro-scale
imaging of biological tissue is based on quantitative phase retrieval
schemes. By incorporating refraction, this method yields improved contrast
compared to purely absorption-based radiography but involves a phase
retrieval problem since of physical limitation of detectors. A general
method is proposed in this paper for one step reconstruction of the ray
integral of complex refractive index of an optically weak object from
intensity distribution of the hologram.


\textbf{AMS }subject classification\textbf{\ }78A45 78A46

\section{Introduction}

A sample is illuminated by a parallel beam of coherent X-rays; the intensity
of the diffracted pattern is registered at a distant plane detector for
determination of distribution of attenuation and refractivity in the sample.
By incorporating refraction, this method yields improved contrast compared
to purely absorption-based radiography but involves a phase retrieval
problem. The linearized version of\ this problem is known as the contrast
transfer function model (CTF). This model is applied to optically weak
objects and the Helmholtz equation is replaced by the paraxial approximation 
\cite{Cl},\cite{PMGMW}. In \cite{JL} the near-field phase retrieval problem
was considered for general compactly supported objects, given two
independent intensity patterns recorded at different distances or
frequencies. A version of the Newton method was applied for 2D
reconstruction from one hologram \cite{MBKH}. The method of phase contrast
imaging was applied for resolving the refractive index in three dimensions 
\cite{RS}. Methods of complex analysis were applied in \cite{MH} for proving
stability of the inversion. An exact estimate of the norm of the inversion
is given showing exponential growth for large Fresnel numbers \cite{MH}. In
thist paper a theoretically exact method is proposed for reconstruction of
ray integrals of the complex refraction index of object from the Fourier
transform of a single hologram in the frame of CTF model with arbitrary
Fresnel number.

\section{CTF-model}

Intensity of the radiation is the measurable quantity 
\begin{equation*}
I\left( \psi \right) =\left\vert \mathcal{D}\left( \exp \psi \right)
\right\vert ^{2}.
\end{equation*}%
Here $\mathcal{D}$ is the Fresnel propagator (paraxial approximation) in 3D
space generating the near field hologram,. 
\begin{equation*}
\psi =\mathbf{k}\int \left( 1-\mathbf{n}\right) \mathrm{d}z=\mathbf{k}\int
\left( \delta -i\beta \right) \mathrm{d}z,
\end{equation*}%
$\mathbf{n}=1-\delta +i\beta $ is the refraction index of the object, $%
\mathbf{k}$ is the space frequency, and $z$ is the coordinate along the
central ray. For small $\psi ,$ 
\begin{equation*}
\mathcal{D}\left( \exp \psi \right) =1+\mathcal{D}\left( \psi \right)
+O\left( \left\Vert \psi \right\Vert ^{2}\right)
\end{equation*}%
which yields%
\begin{equation*}
I\left( \psi \right) =1+2T\left( \psi \right) +O\left( \left\Vert \psi
\right\Vert \right) \psi
\end{equation*}%
for an appropriate norm\footnote{$2T$ is the operator denoted by $T$ in \cite%
{MH}.
\par
{}}. Operator $T\left( \psi \right) =\func{Re}\mathcal{D}\left( \psi \right) 
$ is the linearisation of $I$ at $\psi =0$ (weak-object approximation).
According to \cite{Pag} the propagator can be written in the form

\begin{equation}
\mathcal{D}(\psi )\doteqdot F^{-1}(m_{\mathfrak{f}}F(\psi )),\ m_{\mathfrak{f%
}}\left( \xi \right) =\exp \left( \frac{-id\left\vert \xi \right\vert ^{2}}{2%
\mathbf{k}}\right) ,  \label{8}
\end{equation}%
where $z=d$ is the distance to detector, and $\left\vert \xi \right\vert
^{2}=\xi _{1}^{2}+\xi _{2}^{2}$ where $\xi =\left( \xi _{1},\xi _{2}\right) $
are coordinates on the frequency plane of the object. We use notation 
\begin{equation*}
\hat{a}\left( \xi \right) =F\left( a\right) =\int \exp \left( -2\pi i\left(
x_{1}\xi _{1}+x_{2}\xi _{2}\right) \right) a\left( x\right) \mathrm{d}x_{1}%
\mathrm{d}x_{2}
\end{equation*}%
for the Fourier transform of a function $a$ on $R^{2}$ where $x=\left(
x_{1},x_{2}\right) $ are coordinates vanishing on the central ray. The
dimensionless Fresnel number is%
\begin{equation}
\mathfrak{f}=\frac{\mathbf{k}b^{2}}{2\pi d},  \label{7}
\end{equation}%
where $b$ is the diameter of a central disc $\Omega _{b}$ in the plane $%
\left\{ z=0\right\} $ that contains the support $S\ $of the projection of
the sample (or of its projection). Integrals%
\begin{equation*}
\mu =\mathbf{k}\int \delta \mathrm{d}z,\ \varphi =\mathbf{k}\int \beta 
\mathrm{d}z
\end{equation*}%
are the phase shift and the attenuation of a ray $\left\{ x=\mathrm{const}%
\right\} $, respectively. The coordinates $y_{i}=x_{i}/b,\ i=1,2$ are
normalized in such way that support of $\psi $ is contained in the disc $%
\left\{ y:\left\vert y\right\vert \leq 1/2\right\} .$ By (\ref{8}) and (\ref%
{7}) we have%
\begin{equation}
\hat{T}\left( \psi \right) =\cos \left( \frac{\pi \left\vert \eta
\right\vert ^{2}}{\mathfrak{f}}\right) \hat{\mu}+\sin \left( \frac{\pi
\left\vert \eta \right\vert ^{2}}{\mathfrak{f}}\right) \hat{\varphi}
\label{1}
\end{equation}%
where $\eta _{i}=b\xi _{i},\ i=1,2$ are coordinates dual to $y.$ Equation (%
\ref{1}) has the unique solution \cite{Mar}. The norm of $T$ is bounded by $1
$ on $L_{2}\left( \mathbb{R}^{2}\right) $ since operator $\mathcal{D}$ is
unitary. The norm of the operator $T^{-1}$ is estimated by $\exp \left( 
\mathrm{const}\ \mathfrak{f}\right) $\ according to \cite{MH}.

\section{Reconstruction of the phase shift}

\textbf{Definition.} A function $Z$ defined on the complex plane is said to
be of \textit{sine-type} if it is entire and

(i) the zero set $\Lambda $ of $Z$ is separated that is there exists $c>0$
such that%
\begin{equation}
\left\vert \lambda -\mu \right\vert \geq c  \label{20}
\end{equation}%
for any $\lambda $,$\mu \in \Lambda ,$

(ii) there are constants $A,B,H$ such that%
\begin{equation}
A\exp \left( \pi \left\vert \func{Im}\lambda \right\vert \right) \geq
\left\vert Z\left( \lambda \right) \right\vert \geq B\exp \left( \pi
\left\vert \func{Im}\lambda \right\vert \right)  \label{4}
\end{equation}%
for any $\lambda $ such that $\left\vert \func{Im}\lambda \right\vert \leq
H. $

We call $Z$ \textit{generating} function for phase shift at a Fresnel number 
$\mathfrak{f}$ if it is of sine-type with real zeros and each zero $\lambda $
satisfies 
\begin{equation}
\lambda ^{2}=\mathfrak{f}\left( l\left( \lambda \right) +\frac{1}{2}\right) 
\label{3}
\end{equation}%
where $l\left( \lambda \right) $ is an integer.

\begin{theorem}
\label{T}If there exists a generating function $Z\ $at a Fresnel number $%
\mathfrak{f\ }$then for arbitrary functions $\varphi ,\mu \in L_{2}\left( 
\mathbb{R}^{2}\right) $ supported by $\Omega _{b},$ the phase shift $\varphi 
$ can be found explicitly from (\ref{1}).
\end{theorem}

\textbf{Remark.} Note that $b$ in (\ref{7}) is arbitrary number that fulfils 
$b\geq b\left( S\right) $ where $b\left( S\right) $ is the minimal diameter
of a central disc that contains $S.$ Therefore the reconstruction can be
done for a sample with an arbitrary Fresnel number $\mathfrak{f}$ just by
choosing closest odd $\mathfrak{f}_{odd}\geq \mathfrak{f}$ and applying
Theorem \ref{T} at the number $\mathfrak{f}_{odd}$ .

\textit{Proof. }Let $P$ denote (Paley-Wiener) space of functions $g\in
L_{2}\left( \mathbb{R}\right) $ such that $F\left( g\right) $ is supported
in $\left[ -1/2,1/2\right] .$

\begin{lemma}
\label{H}Let $Z$ be a generating function at Fresnel number $\mathfrak{f}.$
Then for an arbitrary function $g\in P$ can be interpolated from the set $%
\Lambda $ of roots of $Z$ by%
\begin{equation}
g\left( t\right) =\sum_{\lambda \in \Lambda }\frac{Z\left( t\right) g\left(
\lambda \right) }{\left( t-\lambda \right) Z^{\prime }\left( \lambda \right) 
},  \label{10}
\end{equation}%
where the series converges uniformly on each bounded interval.
\end{lemma}

The proof will be given in the next Section.

\textit{Proof of Theorem. }Let $\theta $ be a unit vector in the frequency
plane. For any zero$\ \lambda ,$ vector $\eta =\pm \lambda \theta $
satisfies $\left\vert \eta \right\vert ^{2}=\mathfrak{f}\left( l\left(
\lambda \right) +1/2\right) $. It follows that $\cos \left( \pi \left\vert
\eta \right\vert ^{2}/\mathfrak{f}\right) =0$ which yields $\hat{\varphi}%
\left( \lambda \theta \right) =\left( -1\right) ^{l\left( \lambda \right) }%
\hat{\Psi}\left( \lambda \theta \right) .$ By the slice theorem for any\
unit vector $\theta ,$%
\begin{equation*}
\hat{\varphi}\left( t\theta \right) =\int \exp \left( -2\pi ipt\right)
R\varphi \left( p,\theta \right) \mathrm{d}p,\ t\in \mathbb{R},
\end{equation*}%
where $R$ means the Radon transform and function $R\varphi \left( \cdot
,\theta \right) \in L_{2}\left( R\right) $ is supported in $\left[ -1/2,1/2%
\right] .$ Therefore $g\in P$ where $g_{\theta }\left( t\right) =\hat{\varphi%
}\left( t\theta \right) .$ Lemma \ref{H} can be applied to $g_{\theta }$ for
any vector $\theta .$ Therefore function $\hat{\varphi}$ is reconstructed by
(\ref{10}) on the whole frequency plane by the formula 
\begin{equation}
\hat{\varphi}\left( \eta \right) =Z\left( \left\vert \eta \right\vert
\right) \sum_{\Lambda }\frac{\left( -1\right) ^{l\left( \lambda \right) }%
\hat{\Psi}\left( \lambda \theta \right) }{\left( \left\vert \eta \right\vert
-\lambda \right) Z^{\prime }\left( \lambda \theta \right) }  \label{16}
\end{equation}%
where $\theta =\left\vert \eta \right\vert ^{-1}\eta $ and$\ \Psi =T\left(
\psi \right) .$ Application of the inverse FT recovers the phase shift $%
\varphi .$ $\square $

\textbf{Remark. }If a generating function $Z$ is even function (\ref{16})
can be written in form%
\begin{equation}
\hat{\varphi}\left( \eta \right) =Z\left( \left\vert \eta \right\vert
\right) \sum_{\lambda \in \Lambda ,\lambda >0}^{\infty }\left( -1\right)
^{l\left( \lambda \right) }\left( \frac{\hat{\Psi}\left( \lambda \theta
\right) }{\left( \left\vert \eta \right\vert -\lambda \right) Z^{\prime
}\left( \lambda \theta \right) }+\frac{\hat{\Psi}\left( -\lambda \theta
\right) }{\left( \left\vert \eta \right\vert +\lambda \right) Z^{\prime
}\left( \lambda \theta \right) }\right)  \label{34}
\end{equation}%
since $Z\left( 0\right) \neq 0$.

\section{Interpolation in Paley-Wiener space}

\begin{lemma}
\cite{Y}If $Z$ is a sine-type function, then for arbitrary $\varepsilon >0$
there exists $\delta >0$ such that for any $l$ that is not in $\varepsilon $
neighborhood of a zero of $Z$ inequality holds%
\begin{equation}
\left\vert Z\left( l+ih\right) \right\vert \geq \delta \exp \left( \pi
\left\vert h\right\vert \right) .  \label{19}
\end{equation}
\end{lemma}

\textit{Proof of Lemma \ref{H}. }Let $\Gamma _{\theta }=\left\{ \func{Im}%
\lambda =h\right\} $. We have 
\begin{equation}
\int_{\Gamma _{h}}\left\vert g\left( \lambda \right) \right\vert ^{2}\mathrm{%
d}\lambda \leq \left\Vert g\right\Vert ^{2}\exp \left( 2\pi \left\vert
h\right\vert \right) .  \label{2}
\end{equation}%
since $g\left( \xi +ih\right) =F\left( \exp \left( 2\pi ih\right)
F^{-1}\left( f\right) \right) .$ By (\ref{4}), the Cauchy inequality and (%
\ref{2}) we get for any real $t$, 
\begin{eqnarray*}
B^{2}\left\vert \int_{\Gamma _{h}}\frac{g\left( \lambda \right) }{\left(
\lambda -t\right) Z\left( \lambda \right) }\mathrm{d}\lambda \right\vert
^{2} &\leq &\exp \left( -2\pi \left\vert h\right\vert \right) \left(
\int_{\Gamma _{h}}\left\vert \frac{g\left( \lambda \right) }{\lambda -t}%
\right\vert \mathrm{d}\lambda \right) ^{2} \\
&\leq &\exp \left( -2\pi \left\vert h\right\vert \right) \int_{\Gamma _{h}}%
\frac{\mathrm{d}\lambda }{\left\vert \lambda -t\right\vert ^{2}}\int_{\Gamma
_{h}}\left\vert g\left( \lambda \right) \right\vert ^{2}\mathrm{d}\lambda
\leq \frac{\pi }{\left\vert h\right\vert }\left\Vert g\right\Vert ^{2}.
\end{eqnarray*}%
that is%
\begin{equation*}
\left\vert \int_{\Gamma _{h}}\frac{g\left( \lambda \right) }{\left( \lambda
-t\right) Z\left( \lambda \right) }\mathrm{d}\lambda \right\vert \leq \frac{1%
}{B}\left( \frac{\pi }{h}\right) ^{1/2}\left\Vert g\right\Vert .
\end{equation*}%
For any $l>0,$ one can replace here $\Gamma _{h}$ by $\Gamma _{h,l}=\left\{
\left\vert \func{Re}\lambda \right\vert \leq l,\ \left\vert \func{Im}\lambda
\right\vert =h,\right\} $. This implies

\begin{equation*}
\left\vert \int_{\Gamma _{h,l}}\frac{g\left( \lambda \right) }{\left(
\lambda -t\right) Z\left( \lambda \right) }\mathrm{d}\lambda \right\vert
\rightarrow 0
\end{equation*}%
as $h\rightarrow \infty $ uniformly with respect to $l.$\ Set $\Gamma
_{l,h}=\left\{ \func{Re}\lambda =\pm l,\ \left\vert \func{Im}\lambda
\right\vert \leq h\right\} $ and choose a sequence $l=l_{m}\rightarrow
\infty $ such that $l_{m}\ $and $-l_{m}$ are kept at a distance $\geq c/4$
from $\Lambda $ for any $m$ where $c$ is the constant in (\ref{20}). We have
then%
\begin{equation*}
\left\vert \int_{\Gamma _{l,h}}\frac{g\left( \lambda \right) }{\left(
\lambda -t\right) Z\left( \lambda \right) }\mathrm{d}\lambda \right\vert
\rightarrow 0
\end{equation*}%
as $k\rightarrow \infty $ uniformly for $h.$ This property follows from (\ref%
{19}) if we take $\varepsilon =c/4.$ The sum%
\begin{equation}
\int_{\Gamma _{l,h}}+\int_{\Gamma _{h,l}}\frac{g\left( \lambda \right) }{%
\left( \lambda -t\right) Z\left( \lambda \right) }\mathrm{d}\lambda
\label{21}
\end{equation}%
is the integral over the perimeter of the rectangle of size $2l\times 2h.$
If $l_{k}\rightarrow \infty $ and $h\rightarrow \infty ,$ the sum (\ref{21})
tends to zero. If the orientation of the perimeter is counterclockwise the
sum can be calculated by the Residue theorem for the poles $\lambda =t,\
\left\vert t\right\vert <l,$ $\lambda \in \Lambda $, $\left\vert \lambda
\right\vert <l_{m}.$ This yields%
\begin{equation*}
\mathrm{res}_{t}+\sum \mathrm{res}_{\lambda _{k}}=\frac{g\left( t\right) }{%
Z\left( t\right) }+\sum_{\lambda \in \Lambda ,\left\vert \lambda \right\vert
<l}\frac{g\left( \lambda \right) }{\left( \lambda -t\right) Z^{\prime
}\left( \lambda \right) }\rightarrow 0
\end{equation*}%
as $l\rightarrow \infty .$ The limit of the partial sum vanishes and (\ref%
{10}) follows. $\square $

\section{Small Fresnel numbers}

Here generating functions at $f=1,2,3,4,5$ are constructed.

\textbf{Case} $\mathfrak{f}=1$

Check that $Z_{1}\left( \lambda \right) =\cos \left( \pi \sqrt{\lambda
^{2}-1/4}\right) $ is a generating function. The numbers $\lambda _{k}=\pm 
\sqrt{k^{2}+k+1/2},\ k=0,1,...$ are all zeros and $\lambda _{k}^{2}=l\left(
k\right) +1/2$ where $l\left( k\right) =k^{2}+k.$ By Lemma \ref{H} we obtain
formula%
\begin{equation*}
\hat{\varphi}\left( \eta \right) =\cos \left( \pi \sqrt{\left\vert \eta
\right\vert ^{2}-\frac{1}{4}}\right) \sum_{0}^{\infty }\left( -1\right)
^{l\left( k\right) +k}\frac{k+1/2}{\pi \lambda _{k}}\left( \frac{\hat{\Psi}%
\left( \lambda _{k}\theta \right) }{\left\vert \eta \right\vert -\lambda _{k}%
}+\frac{\hat{\Psi}\left( -\lambda _{k}\theta \right) }{\left\vert \eta
\right\vert +\lambda _{k}}\right)
\end{equation*}%
since $Z_{1}^{\prime }\left( \lambda _{k}\right) =\left( -1\right) ^{k}\pi
\lambda _{k}\left( k+1/2\right) ^{-1}.\ $

\textbf{Case} $\mathfrak{f}=2$

The function $Z_{2}\left( \lambda \right) =\cos \pi \sqrt{\lambda ^{2}-3/4}$
is of the sine-type. The roots of $Z$ fulfil $\lambda _{k}=\left(
k+1/2\right) ^{2}+3/4=k^{2}+k+1=2\left( l\left( k\right) +1/2\right) \ $for $%
k=0,1,2,...$ where $l\left( k\right) =\left( k^{2}+k\right) /2$ is integer%
\textit{\ }for any $k.$\textit{\ }This implies (\ref{3}) at $\mathfrak{f}=2$
hence $Z_{2}$ is a generating function hence Lemma \ref{H} can be applied.

\textbf{Case} $\mathfrak{f}=3$

The function%
\begin{equation}
Z_{0}\left( \lambda \right) =\cos \left( \frac{\pi }{3}\sqrt{\lambda ^{2}+%
\frac{3}{4}}\right) \cos \frac{\pi }{3}\left( \sqrt{\lambda ^{2}-\frac{5}{4}}%
+1\right) \cos \frac{\pi }{3}\left( -\sqrt{\lambda ^{2}-\frac{5}{4}}+2\right)
\label{32}
\end{equation}%
is single-valued and of sine-type. Its zeros are real and equal 
\begin{align*}
\alpha _{k}^{2}& =9\left( k+\frac{1}{2}\right) ^{2}-\frac{3}{4}=3\left(
l\left( k\right) +\frac{1}{2}\right) ,\ \ l\left( k\right) =3k^{2}+3k, \\
\beta _{k}^{2}& =\left( 3\left( k+\frac{1}{2}\right) -1\right) ^{2}+\frac{5}{%
4}=3\left( l\left( k\right) +\frac{1}{2}\right) ,\ \ l\left( k\right)
=3k^{2}+k, \\
\gamma _{k}^{2}& =\left( 3\left( k+\frac{1}{2}\right) -2\right) ^{2}+\frac{5%
}{4}=3\left( l\left( k\right) +\frac{1}{2}\right) ,\ \ l\left( k\right)
=3k^{2}-k,
\end{align*}%
where $k=0,1,2,...$ for the cosine factors in (\ref{32}), respectively. They
satisfy (\ref{3}) at $\mathfrak{f}=3$ however the zeros $\lambda =\pm \sqrt{%
3/2}$ have multiplicity 2. We set%
\begin{equation*}
Z_{3}\left( \lambda \right) =Z_{0}\left( \lambda \right) R\left( \lambda
\right) ,\ R\left( \lambda \right) =\frac{\lambda ^{2}-9/2}{\lambda ^{2}-3/2}%
.
\end{equation*}%
Function $Z_{3}$ is still of sine-type since $R\left( \lambda \right)
\rightarrow 1$ as $\left\vert \func{Im}\lambda \right\vert \rightarrow
\infty ,$ the denominator reduces multiplicity to 1. The numerator has zero $%
\lambda ^{2}=9/2$ which also satisfies (\ref{3}). It follows that $Z$ is a
generating function at $\mathfrak{f}=3.$

\textbf{Case} $\mathfrak{f}=4$

Function%
\begin{equation*}
Z_{0}\left( \lambda \right) =\sin \left( \frac{\pi }{2}\sqrt{\lambda ^{2}-2}%
\right) \cos \frac{\pi }{2}\sqrt{\lambda ^{2}-1}
\end{equation*}%
is of sine-type whose zeros are%
\begin{equation*}
\alpha _{k}=\pm \sqrt{4k^{2}+2};\ \beta _{k}=\pm \sqrt{4\left(
k^{2}+k\right) +2},\ k=0,1,2,....
\end{equation*}%
The zeros fulfil (\ref{3}) and are simple except the double zeros $\alpha
_{0}=\beta _{0}=\pm \sqrt{2}.$ Therefore the product%
\begin{equation*}
Z_{4}=Z_{0}R,\ R\left( \lambda \right) =\frac{\lambda ^{2}-14}{\lambda ^{2}-2%
}
\end{equation*}%
is still of sine-type since $R\left( \lambda \right) \rightarrow 1$ as $%
\left\vert \func{Im}\lambda \right\vert \rightarrow \infty .$ The
denominator reduces multiplicity of the double zeros to 1. The product is a
generating function since all zeros are real and fulfil (\ref{3}). We call $%
R $ \textit{correction} factor.

\textbf{Case} $\mathfrak{f}=5$

Take%
\begin{align}
Z_{5}\left( \lambda \right) & =\cos \frac{\pi }{5}\sqrt{\lambda ^{2}+\frac{15%
}{4}}\cos \frac{\pi }{5}\left( \sqrt{\lambda ^{2}-\frac{1}{4}}+1\right) \cos 
\frac{\pi }{5}\left( -\sqrt{\lambda ^{2}-\frac{1}{4}}+4\right)  \label{26} \\
& \times \cos \frac{\pi }{5}\left( \sqrt{\lambda ^{2}-\frac{9}{4}}+2\right)
\cos \frac{\pi }{5}\left( -\sqrt{\lambda ^{2}-\frac{9}{4}}+3\right) R\left(
\lambda \right) .  \notag
\end{align}%
The zeros of the first five factors are%
\begin{align*}
\alpha _{k}^{2}& =\left( 5\left( k+\frac{1}{2}\right) \right) ^{2}-\frac{15}{%
4}=5\left( l+\frac{1}{2}\right) ,\ l=5k^{2}+5k, \\
\beta _{k}^{2}& =\left( 5\left( k+\frac{1}{2}\right) -1\right) ^{2}+\frac{1}{%
4}=5\left( l+\frac{1}{2}\right) ,\ l=5k^{2}+3k,\  \\
\gamma _{k}^{2}& =\left( 5\left( k+\frac{1}{2}\right) -2\right) ^{2}+\frac{9%
}{4}=5\left( l+\frac{1}{2}\right) ,\ l=5k^{2}+k,\  \\
\delta _{k}^{2}& =\left( 5\left( k+\frac{1}{2}\right) -3\right) ^{2}+\frac{9%
}{4}=5\left( l+\frac{1}{2}\right) ,\ l=5k^{2}-k, \\
\varepsilon _{k}^{2}& =\left( 5\left( k+\frac{1}{2}\right) -4\right) ^{2}+%
\frac{1}{4}=5\left( l+\frac{1}{2}\right) ,\ l=5k^{2}-3k,
\end{align*}%
where $k=0,1,2,...$ and the points $\lambda =\pm \sqrt{5/2}$ $(k=0)$ are
zeros of multiplicity 3. We define the correction factor by%
\begin{equation*}
R\left( \lambda \right) =\frac{\left( \lambda ^{2}-15/2\right) \left(
\lambda ^{2}-35/2\right) }{\left( \lambda ^{2}-5/2\right) ^{2}}.
\end{equation*}%
This make the product (\ref{26}) a generating function at $\mathfrak{f}=5.$

\section{Generating function at arbitrary odd $\mathfrak{f}$}

For an arbitrary odd $\mathfrak{f}=2p+1,\ p>0,$ we take%
\begin{equation*}
Z_{0}\left( \lambda \right) =\cos \left( \frac{\pi }{\mathfrak{f}}\sqrt{\rho
_{0}}\right) \prod_{q=1}^{p}\cos \left( \frac{\pi }{\mathfrak{f}}\left( 
\sqrt{\rho _{q}}+q\right) \right) \cos \left( \frac{\pi }{\mathfrak{f}}%
\left( \sqrt{\rho _{f-q}}+\mathfrak{f}-q\right) \right) ,
\end{equation*}%
where 
\begin{equation*}
\rho _{q}\left( \lambda \right) =\rho _{\mathfrak{f}-q}\left( \lambda
\right) =\lambda ^{2}-\frac{\mathfrak{f}}{2}+\left( \frac{\mathfrak{f}}{2}%
-q\right) ^{2}.
\end{equation*}%
The first factor is a single-valued entire function. The same is true for
the products%
\begin{equation*}
\cos \frac{\pi }{\mathfrak{f}}\left( \sqrt{\rho _{q}}+q\right) \cos \frac{%
\pi }{\mathfrak{f}}\left( \rho _{f-q}+\mathfrak{f}-q\right) ,\ q=1,...,p
\end{equation*}%
since%
\begin{equation*}
\cos \frac{\pi }{\mathfrak{f}}\left( -\sqrt{\rho _{q}}+q\right) =-\cos \frac{%
\pi }{\mathfrak{f}}\left( \sqrt{\rho _{f-q}}+\mathfrak{f}-q\right) .
\end{equation*}%
Therefore the product of all $\mathfrak{f}$ factors is also an single-valued
entire function. Its zeros are%
\begin{equation*}
\lambda _{k,q}^{2}=\mathfrak{f}\left( l+\frac{1}{2}\right) ,\ l=\mathfrak{f}%
k^{2}+\left( \mathfrak{f}-2q\right) k,\ q=0,...,p,\ k=0,1,2,...
\end{equation*}%
where the zero $\lambda =\pm \sqrt{\mathfrak{f/}2}$ appears $p+1$ times. We
define 
\begin{equation}
Z_{\mathfrak{f}}=Z_{0}R.  \label{33}
\end{equation}%
The correction factor%
\begin{equation*}
R\left( \lambda \right) =\frac{\prod_{q=1}^{p}\lambda ^{2}-\mathfrak{f}%
\left( 2q+1/2\right) }{\left( \lambda ^{2}-\mathfrak{f}/2\right) ^{p}}
\end{equation*}%
has zeros $\mu _{q}=\pm \sqrt{\mathfrak{f}\left( 2q+1/2\right) },$ $%
q=1,2,...,p$ that also fulfil (\ref{3}). Check that no of these zeros is a
zero of $Z_{0}.$ Suppose the opposite, let $\lambda _{k,r}^{2}=\mathfrak{f}%
\left( 2q+1/2\right) $ for some $1\leq r,q\leq p.$ Then%
\begin{equation*}
\mathfrak{f}k^{2}+\left( \mathfrak{f}-2r\right) k=2q
\end{equation*}%
for some $k.$ It is not the case if $k=0$ since $q>0.$ In the case $k\geq 1$
we have%
\begin{equation*}
\mathfrak{f}k^{2}=2q-\left( \mathfrak{f}-2r\right) k\leq 2qk-\left( 
\mathfrak{f}-2r\right) k=\left( 2\left( q+r\right) -\mathfrak{f}\right) k<%
\mathfrak{f}k
\end{equation*}%
since $2\left( q+r\right) \leq 4p<2\mathfrak{f.}$ This yields $k^{2}<k$
which is the contradiction. It follows that $Z_{\mathfrak{f}}$ is an even
generating function at $\mathfrak{f}$ for phase shift.\ The zeros of the
main series with $k>0$ fulfil%
\begin{equation}
..,\lambda _{k,p}<\lambda _{k,p-1}<...<\lambda _{k,0}<\lambda
_{k+1,p}<\lambda _{k+1,p}<...  \label{5}
\end{equation}%
and 
\begin{equation}
\lambda _{k,q}=\mathfrak{f}\left( k+\frac{1}{2}\right) -q+O\left(
k^{-1}\right) .  \label{9}
\end{equation}%
It follows that the gap between adjacent zeros tends to 1 as $k\rightarrow
\infty .$

Formulas for the cases $\mathfrak{f}=1,3,5$ as above are particular cases of
(\ref{33}).

\begin{conclusion}
Interpolation formula for the phase shift like (\ref{34}) holds for $Z_{%
\mathfrak{f}}$ for any odd $\mathfrak{f}$.
\end{conclusion}

\textbf{Remark. }A similar construction can be applied for even Fresnel
numbers.

\section{Reconstruction of the attenuation}

The similar method can be applies to$\ $reconstruction of the attenuation $%
\mu .$ We call a function $W$ generating at a Fresnel number $\mathfrak{f}$
for attenuation. if it is of sine-type, all zeros $\lambda $ are real and
fulfil the equation%
\begin{equation}
\left\vert \lambda \right\vert ^{2}=\mathfrak{f}l\left( \lambda \right)
\label{30}
\end{equation}%
where $l$ is an integer for any $\lambda .$ Function%
\begin{equation*}
W_{1}\left( \lambda \right) =\cos \left( \pi \sqrt{\lambda ^{2}+\frac{1}{4}}%
\right)
\end{equation*}%
is generating at Fresnel number $\mathfrak{f}=1$ for the attenuation with\
zeros $\lambda _{k}=\pm \sqrt{k^{2}+k},$ $k=0,1,2,....\ $A generating
function for attenuation can be constructed at any odd Fresnel number $%
\mathfrak{f}$ as follows. We can take 
\begin{equation*}
W_{\mathfrak{f}}\left( \lambda \right) =\cos \frac{\pi }{\mathfrak{f}}\sqrt{%
\sigma _{0}}\prod_{q=1}^{p}\cos \frac{\pi }{\mathfrak{f}}\left( \sqrt{\sigma
_{q}}+q\right) \cos \frac{\pi }{\mathfrak{f}}\left( \sqrt{\sigma _{\mathfrak{%
f-}q}}+\mathfrak{f}-q\right) R\left( \lambda \right)
\end{equation*}%
where%
\begin{equation*}
\sigma _{q}\left( \lambda \right) =\lambda ^{2}+\left( \frac{\mathfrak{f}}{2}%
-q\right) ^{2},\ q=0,...,\mathfrak{f}-1.
\end{equation*}%
The zeros of the first $\mathfrak{f}$ factors are%
\begin{equation}
\lambda _{k,q}^{2}=\mathfrak{f}l,\ l=\mathfrak{f}k^{2}+\left( \mathfrak{f}%
-2q\right) k,\ q=0,...,p,\ k=0,1,2,...  \label{31}
\end{equation}%
They fulfil (\ref{30}) but number $\lambda =0$ appears $p+1$ times. The
factor%
\begin{equation*}
R\left( \lambda \right) =\frac{\prod_{q=0}^{p-1}\left( \lambda ^{2}-\left(
2q+1\right) \mathfrak{f}\right) }{\lambda ^{2p}}
\end{equation*}%
corrects to 1 the multiplicity of this zero. The zeros of the numerator
satisfies (\ref{30}) and do not coincide with points (\ref{31}). Therefore
function \thinspace $W_{\mathfrak{f}}$ generates attenuation and the
arguments of Theorem \ref{T} work for reconstruction of attenuation $\mu $.
Note that points (\ref{31}) satisfies inequalities (\ref{5}) and (\ref{9}).

\begin{conclusion}
Interpolation formula for attenuation like (\ref{34}) holds for $W_{%
\mathfrak{f}}$ for any odd $\mathfrak{f}$.
\end{conclusion}

\section{Convergence and truncation error}

\begin{theorem}
If $\Lambda $ is the set of zeros of a function $Z$ of sine-type, then the
functions%
\begin{equation*}
\frac{Z\left( t\right) }{\left( t-\lambda \right) Z^{\prime }\left( \lambda
\right) },\ \lambda \in \Lambda
\end{equation*}%
form a\ Riesz basis in the space $P.$ This means that $\left\{ g\left(
\lambda \right) \right\} \in l_{2}\ $for any function $g\in P$ and the
sequence(\ref{10}) converges to $g$ in $L_{2}\left( \mathbb{R}\right) .$
Vice versa\ for any sequence $\left\{ c_{\lambda }\right\} \in l_{2}$ the
series%
\begin{equation}
\sum_{\lambda \in \Lambda }g\left( \lambda \right) \frac{Z\left( t\right) }{%
\left( t-\lambda \right) Z^{\prime }\left( \lambda \right) }  \label{29}
\end{equation}%
converges to an element of $P.$
\end{theorem}

Pointwise convergence of series (\ref{29}) can be evaluated for the
classical Whittecker-Kotelnikov-Shannon series which is the particular case
for $Z\left( t\right) =\sin \pi t.$ Take the indicator function $f$ of the
set $\left[ -1,-2/3\right] \cup \left[ 2/3,1\right] $ as a model for the
space of piecewise smooth functions. Compare function $g=\hat{f}\ $of
Paley-Wiener type with the truncated WKS sum 
\begin{equation*}
g_{N}\left( t\right) =\sum_{k=-N}^{N}g\left( k\right) \frac{\sin \left( \pi
\left( t-k\right) \right) }{\pi \left( t-k\right) }.
\end{equation*}%
A calculation shows that the partial sum $g_{8}$ gives rather good
approximation: 
\begin{equation*}
\max_{\left\vert t\right\vert \leq 6}\left\vert g\left( t\right)
-g_{8}\left( t\right) \right\vert \leq 0.006.
\end{equation*}%
The series (\ref{29}) is of the same as the WKS series. Therefore the
convergence of these series is expected to have the same rate.


\begin{thebibliography}{9}
\bibitem{Cl} P. Cloetens, W. Ludwig, J. Baruchel, D. Van Dyck, J. Van
Landruyt, J. P. Guigay, and M. Schlemker, \textit{Holotomography: Quantative
phase tomography with micrometer resolution using hard synchrotron radiation
X rays, }Appl. Phys. Lett. \textbf{75 }(1999), 2912-2814.

\bibitem{JL} P. Jonas and A. Louis, \textit{Phase contrast tomography using
holographic measurements,} Inverse Probl.\textit{\ }\textbf{20 }(2014), 75.

\bibitem{Mar} S. Maretzke, \textit{A uniqueness result for propagation-based
phase contrast imaging from a single measurement,} Inverse Probl.\textit{\ }%
\textbf{31} (2015), 065003.

\bibitem{MBKH} S. Maretzke, N. Bartels, M. Krenkel, T. Salditt and T.
Hohage, \textit{Regularized Newton methods for X-ray phase contrast and
general imaging problems,} Opt. Express \textbf{24} (2016), 6490-6506.

\bibitem{MH} S. Maretzke and T. Hohage, \textit{Stability estimates for
linearized near-field phase retrieval in X-ray phase contrast imaging,} SIAM
J. Appl. Math. \textbf{77} (2017), 384-408.

\bibitem{PMGMW} D. Paganin, S. Mayo, T. E. Gureyev, P. R. Miller and \ S. W.
Wilkins, \textit{Simultaneous phase and amplitude extraction from a simple
defocused image of a homogeneous object, }J. Microsc.\textit{\ }\textbf{206 }%
(2002), 33-40.

\bibitem{Pag} D. Paganin, \textit{Coherent X-Ray Optics}, Vol. 1, Oxford
University Press, Oxford, UK, 2006.

\bibitem{RS} A. Ruhlandt and T. Salditt, \textit{Three-dimensional
propagation in near-field tomographic X-ray retrieval,} Acta Crystallogr.
Sect. A \textbf{72} (2016), 215-221.

\bibitem{Y} Young R M 1980 \textit{An introduction to nonharmonic Fourier
series,} New York: Academic Press.
\end{thebibliography}
\end{document}